\documentclass[a4paper,11pt]{amsart}
\usepackage{amssymb}

\input xypic
\xyoption{all}
\xyoption{arc}

\newtheorem{theo}{Theorem}[section]

\newtheorem{prop}[theo]{Proposition}

\newtheorem{coro}[theo]{Corollary}

\def\remark#1{{\refstepcounter{theo}\label{#1}\noindent\sc Remark  
\arabic{section}.\arabic{theo} - }}
\def\example#1{{\refstepcounter{theo}\label{#1}\noindent\sc Example 
\arabic{section}.\arabic{theo} - }}

\def\equat{\refstepcounter{theo}$$~}
\def\endequat{\leqno{\boldsymbol{(\arabic{section}.\arabic{theo})}}~$$}

\def\AG{{\mathfrak A}}
\def\CM{{\mathbb{C}}}
\def\FM{{\mathbb{F}}}

\def\QM{{\mathbb{Q}}}
\def\pG{{\mathfrak p}}
\def\SG{{\mathfrak S}}
\def\ZM{{\mathbb{Z}}}
\def\RM{{\mathbb{R}}}
\def\OC{{\mathcal{O}}}
\def\RC{{\mathcal{R}}}
\def\DC{{\mathcal{D}}}
\def\SC{{\mathcal{S}}}
\def\PC{{\mathcal{P}}}

\def\g{\gamma}

\def\d{\delta}

\def\e{\varepsilon}
\def\ph{\varphi}
\def\l{\lambda}

\def\o{\omega}

\def\th{\theta}

\def\t{\tau}
\def\z{\zeta}
\def\eba{{\bar{e}}}
\def\tba{{\bar{t}}}
\def\fba{{\bar{f}}}
\def\fonction#1#2#3#4#5{\begin{array}{rccc}
{#1} : & {#2} & \longto & {#3} \\
& {#4} & \longmapsto & {#5} 
\end{array}}

\DeclareMathOperator{\Class}{{\mathrm{Class}}}
\DeclareMathOperator{\Ind}{{\mathrm{Ind}}}
\DeclareMathOperator{\Irr}{{\mathrm{Irr}}}
\DeclareMathOperator{\Ker}{{\mathrm{Ker}}}

\DeclareMathOperator{\Res}{{\mathrm{Res}}}
\DeclareMathOperator{\Rad}{{\mathrm{Rad}}}
\DeclareMathOperator{\Id}{{\mathrm{Id}}}
\DeclareMathOperator{\res}{{\mathrm{res}}}
\DeclareMathOperator{\ind}{{\mathrm{ind}}}
\def\epsba{{\bar{\e}}}
\def\chiba{{\bar{\chi}}}
\def\to{\rightarrow}
\def\longto{\longrightarrow}

\def\DS{\displaystyle}
\def\SS{\scriptstyle}

\def\finl{~$\SS \square$}

\def\lexp#1#2{\kern\scriptspace\vphantom{#2}^{#1}\kern-\scriptspace#2}
\def\le{\hspace{0.1em}\mathop{\leqslant}\nolimits\hspace{0.1em}}
\def\ge{\hspace{0.1em}\mathop{\geqslant}\nolimits\hspace{0.1em}}

\def\eqna{\begin{eqnarray*}}
\def\endeqna{\end{eqnarray*}}

\def\itemth#1{\item[${\mathrm{(#1)}}$]}

\def\loewy{\ell}
\def\eval{{\mathrm{ev}}}
\def\evalb{{\overline{\mathrm{ev}}}}

\DeclareMathOperator{\degbar}{{\overline{\mathrm{deg}}}}
\DeclareMathOperator{\ext}{{\mathrm{ext}}}
\DeclareMathOperator{\Ext}{{\mathrm{Ext}}}
\DeclareMathOperator{\Part}{{\mathrm{Part}}}
\DeclareMathOperator{\Log}{{\mathrm{Log}}}

\DeclareMathOperator{\proj}{{\mathrm{proj}}}
\DeclareMathOperator{\projbar}{{\overline{\mathrm{proj}}}}
\DeclareMathOperator{\Resbar}{{\overline{\mathrm{Res}}}}

\def\vertical{\vphantom{\DS{A_A^A}}}
\def\ssim{\!\!\sim}
\def\groth{\RC}
\def\gfp{{\FM_{\! p}}}

\begin{document}

\baselineskip16pt

\title{On the character ring \\
of a finite group}

\author{C\'edric Bonnaf\'e}
\address{
Laboratoire de Math\'ematiques de Besan\c{c}on (CNRS: UMR 6623), 
Universit\'e de Franche-Comt\'e, 
16 Route de Gray, 
25030 Besan\c{c}on Cedex, 
France}

\makeatletter
\email{bonnafe@math.univ-fcomte.fr}
\urladdr{http://www-math.univ-fcomte.fr/pp\_Annu/CBONNAFE/}

\makeatother


\date{\today}

\begin{abstract} 
Let $G$ be a finite group and let $k$ be a sufficiently large 
finite field. Let $\groth(G)$ denote the character ring of $G$ 
(i.e. the Grothendieck ring of the category of $\CM G$-modules). We study 
the structure and the representations of the commutative algebra 
$k \otimes_\ZM \groth(G)$. 
\end{abstract}

\maketitle

\pagestyle{myheadings}

\markboth{\sc C. Bonnaf\'e}{\sc Character ring of a finite group}

\bigskip


Let $G$ be a finite group. We denote by $\groth(G)$ the {\it Grothendieck ring} 
of the category of $\CM G$-modules (it is usually called the {\it character 
ring} of $G$). It is a natural question to try to recover properties 
of $G$ from the knowledge of $\groth(G)$. It is clear that two finite 
groups having the same character table have the same Grothendieck rings and 
it is a Theorem of Saksonov \cite{saksonov} that the converse also holds. 
So the problem is reduced to an intensively studied question in character 
theory: recover properties of the group through properties of its character 
table. 

In this paper, we study the $k$-algebra $k\groth(G)=k \otimes_\ZM \groth(G)$, where 
$k$ is a splitting field for $G$ of positive characteristic $p$. It is clear 
that the knowledge of $k\groth(G)$ is a much weaker information than 
the knowledge of $\groth(G)$. The aim of this paper is to gather results 
on the representation theory of the algebra $k\groth(G)$: although most of 
the results are certainy well-known, we have not found any general treatment 
of these questions. The blocks of $k\groth(G)$ are local algebras which are 
parametrized by conjugacy classes of $p$-regular elements of $G$. 
So the simple $k\groth(G)$-modules are parametrized by 
conjugacy classes of $p$-regular elements of $G$. Moreover, the dimension of the 
projective cover of the simple module associated to the conjugacy class of 
the $p$-regular element $g \in G$ is equal to the number of conjugacy classes 
of $p$-elements in the centralizer $C_G(g)$. We also prove that 
the radical of $k\groth(G)$ is the kernel of the decomposition map 
$k\groth(G) \to k\otimes_\ZM \groth(kG)$, where $\groth(kG)$ is the Grothendieck 
ring of the category of $kG$-modules (i.e. the ring of virtual Brauer characters of $G$). 

We prove that the block of $k\groth(G)$ 
associated to the $p'$-element 
$g$ is isomorphic to the block of $k\groth(C_G(g))$ associated to $1$ 
(such a block is called the {\it principal block}). 
This shows that the study of blocks of $k\groth(G)$ is reduced to the study of 
principal blocks. We also show that the principal block of $k\groth(G)$ 
is isomorphic to the principal block of $k\groth(H)$ whenever $H$ is 
a subgroup of $p'$-index which controls the fusion of $p$-elements 
or whenever $H$ is the quotient of $G$ by a normal $p'$-subgroup. 

We also introduce several numerical invariants (Loewy length, dimension of 
${\mathrm{Ext}}$-groups) that are partly related to the structure of $G$. 
These numerical invariants are computed completely whenever $G$ is the symmetric 
group $\SG_n$ (this relies on previous work of the author: the 
descending Loewy series of $k\groth(\SG_n)$ was entirely computed in \cite{bonnafe sym}) 
or $G$ is a dihedral group and $p=2$. We also provide tables for these 
invariants for small groups (alternating groups $\AG_n$ with $n \le 12$, 
some small simple groups, groups $PSL(2,q)$ with $q$ a prime power $\le 27$, 
exceptional finite Coxeter groups).

\bigskip

\noindent{\sc Notation - } 
Let $\OC$ be a Dedekind domain of characteristic zero, let 
$\pG$ be a maximal ideal of $\OC$, let $K$ be the fraction field 
of $\OC$ and let $k=\OC/\pG$. Let $\OC_\pG$ be the localization of 
$\OC$ at $\pG$: then $k=\OC_\pG/\pG\OC_\pG$. 
If $x \in \OC_\pG$, we denote by $\bar{x}$ its image in $\OC_\pG/\pG \OC_\pG=k$. 
Throughout this paper, we assume that $k$ has characteristic $p > 0$ and 
that $K$ and $k$ are splitting fields for all the finite groups involved 
in this paper. If $n$ is a non-zero natural number, $n_{p'}$ denotes the largest 
divisor of $n$ prime to $p$ and we set $n_p=n/n_{p'}$. 

If $F$ is a field and if $A$ is a finite dimensional 
$F$-algebra, we denote by $\groth(A)$ its Grothendieck group. 
If $M$ is an $A$-module, the radical of $M$ is denoted 
by $\Rad M$ and the class of $M$ in $\groth(A)$ is denoted by $[M]$. 
If $S$ is a simple $A$-module, we denote by $[M:S]$ the multiplicity 
of $S$ as a chief factor of a Jordan-H\"older series of $M$. 
The set of irreducible characters of $A$ is denoted by $\Irr A$.

We fix all along this paper a finite group $G$. 
For simplification, we set $\groth(G)=\groth(KG)$ and 
$\Irr G=\Irr KG$ 
(recall that $K$ is a splitting field for $G$). The abelian group 
$\groth(G)$ is endowed with a structure of ring induced by 
the tensor product. If $\chi \in \groth(G)$, we denote by $\chi^*$ its dual 
(as a class function on $G$, we have $\chi^*(g)=\chi(g^{-1})$ for any $g \in G$). 
If $R$ is any commutative ring, 
we denote by $\Class_R(G)$ the space of class functions $G \to R$ 
and we set $R\groth(G)=R \otimes_\ZM \groth(G)$. If $X$ is a subset of $G$, 
we denote by $1_X^R : G \to R$ the characteristic function of $X$. 
If $R$ is a subring of $K$, then we simply write $1_X=1_X^R$. 
Note that $1_G$ is the trivial character of $G$. 
If $f$, $f' \in \Class_K(G)$, we set 
$$\langle f,f'\rangle_G = \frac{1}{|G|}\sum_{g \in G} f(g^{-1})f'(g).$$
Then $\Irr G$ is an orthonormal basis of $\Class_K(G)$. 
We shall identify $\groth(G)$ with the sub-$\ZM$-module (or sub-$\ZM$-algebra) 
of $\Class_K(G)$ generated by $\Irr G$, and $K\groth(G)$ with $\Class_K(G)$. 
If $f \in \OC_\pG\groth(G)$, we denote by $\fba$ its image in 
$k\groth(G)$.

If $g$ and $h$ are two elements of $G$, we write $g \sim h$ 
(or $g \sim_G h$ if we need to emphasize the group) 
if they are conjugate in $G$. We denote by $g_p$ (resp. $g_{p'}$) 
the $p$-part (resp. the $p'$-part) of $g$. If $X$ is a subset 
of $G$, we set $X_{p'}=\{g_{p'}~|~g \in X\}$ and $X_p=\{g_p~|~g \in X\}$. 
If moreover $X$ is closed under conjugacy, the set of conjugacy classes 
contained in $X$ is denoted by $X/\ssim$. In this case, 
$1_X^R \in \Class_R(G)$. The centre of $G$ is denoted by $Z(G)$. 

\bigskip

\section{Preliminaries}

\bigskip

\subsection{Symmetrizing form} 
Let 
$$\fonction{\t_G}{\groth(G)}{\ZM}{\chi}{\langle \chi, 1_G\rangle_G}$$
denote the canonical symmetrizing form on $\groth(G)$. The dual 
basis of $\Irr G$ is $(\chi^*)_{\chi \in \Irr G}$. It is then 
readily seen that $(\groth(G),\Irr G)$ is a {\it based ring} 
(in the sense of Lusztig \cite[Page 236]{lusztig}). 

If $R$ is any ring, we denote by $\t_G^R : R\groth(G) \to R$ 
the symmetrizing form $\Id_R \otimes_\ZM \t_G$. 

\bigskip

\subsection{Translation by the centre} 
If $\chi \in \Irr G$, we denote by $\o_\chi : Z(G) \to \OC^\times$ the 
linear character such that $\chi(zg)=\o_\chi(z) \chi(g)$ for all 
$z \in Z(G)$ and $g \in G$. If $z \in Z(G)$, we denote by 
$t_z : K\groth(G) \to K\groth(G)$ the linear map defined by $(t_z f)(g)=f(zg)$ 
for all $f \in K\groth(G)$ and $g \in G$. It is clear that 
$t_{zz'}=t_z \circ t_{z'}$ for all $z$, $z' \in Z(G)$ and that $t_z$ 
is an automorphism of algebra. Moreover, 
$$t_z \chi = \o_\chi(z) \chi$$ 
for every $\chi \in \Irr G$. Therefore, $t_z$ is an isometry 
which stabilizes $\OC\groth(G)$. If $R$ is a subring of $K$ such that 
$\OC \subset R \subset K$, we still denote by 
$t_z : R\groth(G) \to R\groth(G)$ the restriction of $t_z$. 
Let $\tba_z = \Id_k \otimes_\OC t_z : k\groth(G) \to k\groth(G)$. This 
is again an automorphism of $k$-algebra. 
If $z$ is a $p$-element, then $\tba_z=\Id_{k\groth(G)}$. 

\bigskip

\subsection{Restriction} 
If $\pi : H \to G$ is a morphism 
of groups, then the {\it restriction through $\pi$} induces 
a morphism of rings $\Res_\pi : \groth(G) \to \groth(H)$. 
If $R$ is a subring of $K$, we still denote by 
$\Res_\pi : R\groth(G) \to R\groth(H)$ the morphism $\Id_R \otimes_\ZM \Res_\pi$. 
We denote by $\Resbar_\pi : k\groth(G) \to k\groth(H)$ 
the reduction modulo $\pG$ of $\Res_\pi : \OC\groth(G) \to \OC\groth(H)$. 
Recall that, if $H$ is a subgroup of $G$ and $\pi$ is the canonical 
injection, then $\Res_\pi$ is just $\Res_H^G$. In this case, 
$\Resbar_\pi$ will be denoted by $\Resbar_H^G$. Note the following fact:
\equat\label{injectif surjectif}
\text{\it If $\pi$ is surjective, then $\Resbar_\pi$ is injective.}
\endequat
\begin{proof}[Proof of \ref{injectif surjectif}]
Indeed, if $\pi$ is surjective, then $\Res_\pi : \groth(G) \to \groth(H)$ is 
injective and its image is a direct summand of $\groth(H)$. 
\end{proof}

\bigskip

\subsection{Radical} 
First, note that, since $k\groth(G)$ is commutative, we have
\equat\label{nilpotent}
\text{\it $\Rad k\groth(G)$ is the ideal of nilpotent elements of $k\groth(G)$.}
\endequat
So, if $\pi : H \to G$ is a morphism of finite groups, then 
\equat\label{stable}
\Resbar_\pi\bigl(\Rad k \groth(G)\bigr) \subset \Rad k\groth(H).
\endequat
The {\it Loewy length} of the algebra $k\groth(G)$ is defined as the smallest natural number 
$n$ such that $\bigl(\Rad k\groth(G)\bigr)^n=0$. We denote it by $\loewy_p(G)$. 
By \ref{injectif surjectif} and \ref{stable}, we have:
\equat\label{surjectif croissant}
\text{\it If $\pi$ is surjective, then $\loewy_p(G) \le \loewy_p(H)$.}
\endequat

\bigskip

\section{Modules for $K\groth(G)$ and $k\groth(G)$}

\medskip

\subsection{Semisimplicity}
Recall that $K\groth(G)$ is identified with the algebra 
of class functions on $G$. 
If $C \in G/\ssim$ and $f \in K\groth(G)$, we denote by $f(C)$ the 
constant value of $f$ on $C$. 
We now define $\eval_C : K\groth(G) \to K$, $f \mapsto f(C)$. 
It is a morphism of $K$-algebras. In other words, it is an 
irreducible representation (or character) of $K\groth(G)$. 
We denote by $\DC_C$ the corresponding simple $K\groth(G)$-module 
($\dim_K \DC_C = 1$ and an element $f \in K\groth(G)$ acts on $\DC_C$ 
by multiplication by $\eval_C(f)=f(C)$). 
Now, $1_C$ is a primitive idempotent of $K\groth(G)$ and 
it is easily checked that 
\equat
K\groth(G) 1_C \simeq \DC_C.
\endequat
Recall that 
\equat\label{formule ec}
1_C=\frac{|C|}{|G|}\sum_{\chi \in \Irr G} \chi(C^{-1}) \chi
\endequat
and
\equat\label{somme ec}
\sum_{C \in G/\!\sim} 1_C = 1_G.
\endequat
Therefore:

\begin{prop}\label{simples 0}
We have:
\begin{itemize}
\itemth{a} $(\DC_C)_{C \in G/\!\sim}$ is a family of representatives 
of isomorphy classes of simple $K\groth(G)$-modules.

\itemth{b} $\Irr K\groth(G)=\{\eval_C~|~C \in G/\ssim\}$.

\itemth{c} $K\groth(G)$ is split semisimple.
\end{itemize}
\end{prop}

We conclude this section by the computation of the Schur elements 
(see \cite[7.2]{geck} for the definition) associated to each irreducible 
character of $K\groth(G)$. Since 
\equat\label{tau GK}
\t_G^K = \sum_{C \in G/\!\sim} \frac{|C|}{|G|}\eval_C,
\endequat
we have by \cite[Theorem 7.2.6]{geck}:

\begin{coro}\label{schur}
Let $C \in G/\ssim$. Then the Schur element associated with the irreducible 
character $\eval_C$ is $\DS{\frac{|G|}{|C|}}$.
\end{coro}

\remark{tz unc}
If $z \in Z(G)$, then $t_z$ induces an isomorphism of algebras 
$K\groth(G)1_C \simeq K\groth(G)1_{z^{-1} C}$.\finl

\bigskip

\remark{expression} 
If $f \in K\groth(G)$, then $f=\sum_{C \in G/\ssim} f(C) 1_C$.\finl

\bigskip

\example{degre} 
The map $\eval_1$ will sometimes be denoted by $\deg$, since it sends 
a character to its degree.\finl

\bigskip

\def\DCba{\bar{\DC}}
\def\PCba{\bar{\PC}}
\def\phba{\bar{\varphi}}

\subsection{Decomposition map} 
Let $d_\pG : \groth(G) \to \groth(kG)$ denote the decomposition map. 
If $R$ is any commutative ring, we denote by $d_\pG^R : R\groth(G) \to R\groth(kG)$ 
the induced map. Note that $\groth(kG)$ is also a ring (for the multiplication 
given by tensor product) and that $d_\pG$ is a morphism of ring. Also, 
by \cite[Corollary 18.14]{curtis}, 
\equat\label{surjectif d}
\text{\it $d_\pG$ is surjective.}
\endequat
Since $\Irr(kG)$ is a linearly independent family of 
class functions $G \to k$ (see \cite[Theorem 17.4]{curtis}), the map 
$\chi : k\groth(kG) \to \Class_k(G)$ that sends the class of a $kG$-module 
to its character is (well-defined and) injective. This is a morphism 
of $k$-algebras. 

Now, if $C$ is a conjugacy class of $p$-regular elements 
(i.e. $C \in G_{p'}/\ssim$), we define 
$$\SC_{p'}(C)=\{g \in G~|~g_{p'}\in C\}$$
(for instance, $\SC_{p'}(1)=G_p$). 
Then $\SC_{p'}(C)$ is called the {\it $p'$-section} of $C$: 
this is a union of conjugacy classes of $G$. 
Let $\Class_k^{p'}(G)$ be the space of class functions $G \to k$ 
which are constant on $p'$-sections. Then, by \cite[Lemma 17.8]{curtis}, 
$\Irr(kG) \subset \Class_k^{p'}(G)$, 
so the image of $\chi$ is contained in $\Class_k^{p'}(G)$. But, 
$\chi$ is injective, $|\Irr(kG)|=|G_{p'}/\!\sim|$ (see \cite[Corollary 17.11]{curtis}) 
and $\dim_k \Class_k^{p'}(G)=|G_{p'}/\ssim|$. Therefore, we can identify, through $\chi$, 
the $k$-algebras $k\groth(kG)$ and $\Class_k^{p'}(G)$. In particular, 
\equat\label{split}
\text{\it $k\groth(kG)$ is split semisimple.}
\endequat

\bigskip

\subsection{Simple ${\boldsymbol{k\groth(G)}}$-modules}
If $C \in G/\ssim$, we still denote by $\eval_C : \OC\groth(G) \to \OC$ the 
restriction of $\eval_C$ and we denote by $\evalb_C : k \groth(G) \to k$ the 
reduction modulo $\pG$ of $\eval_C$. It is easily checked that $\evalb_C$ 
factorizes through the decomposition map $d_\pG$. Indeed, if 
$\eval_C^k : k\groth(kG) \to k$ denote the evaluation at $C$ 
(recall that $k\groth(kG)$ is identified, via the map $\chi$ 
of the previous subsection, to $\Class_k^{p'}(G)$), then 
\equat\label{factorisation}
\evalb_C = \eval_C^k \circ d_\pG^k.
\endequat
Let $\bar{\DC}_C$ be the 
corresponding simple $k\groth(G)$-module. 
Let $\delta_\pG : \groth(K\groth(G)) \to \groth(k\groth(G))$ denote 
the decomposition map (see \cite[7.4]{geck} for the definition). Then
\equat\label{decomposition}
\delta_\pG [\DC_C]=[\bar{\DC}_C].
\endequat
The following facts are well-known:

\begin{prop}\label{isomorphisme}
Let $C$, $C' \in G/\ssim$. Then $\DCba_C \simeq \DCba_{C'}$ if and only if 
$C_{p'} = C_{p'}'$. 
\end{prop}

\begin{proof}
The ``if'' part follows from the following classical fact 
\cite[Proposition 17.5 (ii) and (iv) and Lemma 17.8]{curtis}: if 
$\chi \in \groth(G)$ and if $g \in G$, then 
$$\chi(g) \equiv \chi(g_{p'}) \mod \pG.$$
The ``only if'' part follows from \ref{factorisation} and from the 
surjectivity of the decomposition map $d_\pG$.
\end{proof}

\begin{coro}\label{simples}
We have:
\begin{itemize}
\itemth{a} $(\DCba_C)_{C \in G_{p'}/\!\sim}$ is a family of representatives 
of isomorphy classes of simple $k\groth(G)$-modules.

\itemth{b} $\Irr k\groth(G)=\{\evalb_C~|~C \in G_{p'}/\ssim\}$.

\itemth{c} $\Rad k\groth(G)=\Ker d_\pG^k$.

\itemth{d} $k\groth(G)$ is split.
\end{itemize}
\end{coro}

\begin{proof}
(a) follows from \ref{decomposition} and from the fact that the isomorphy class of 
any simple $k\groth(G)$-modules must occur in some $\d_\pG [S]$, where 
$S$ is a simple $K\groth(G)$-module. (b) follows from (a). 
(c) and (d) follow from (a), (b), \ref{factorisation} and \ref{split}.
\end{proof}

\begin{coro}\label{dimension radical}
$\dim_k \Rad(k\groth(G))=|G/\ssim|-|G_{p'}/\ssim|$.
\end{coro}

\begin{coro}\label{critere semisimple}
$k\groth(G)$ is semisimple if and only if $p$ does not divide $|G|$.
\end{coro}

\bigskip

\example{p group} 
Since $\eval_1$ is also denoted by $\deg$, we shall sometimes denote by 
$\degbar$ the morphism $\evalb_1$. If $G$ is a $p$-group, then 
Corollary \ref{simples} shows that $\Rad k\groth(G) = \Ker(\degbar)$. 
In this case, if $1$, $\l_1$,\dots, $\l_r$ denote the linear characters 
of $G$ and $\chi_1$,\dots, $\chi_s$ denote the non-linear irreducible 
characters of $G$, then $(\overline{\l}_1-1,\dots,\overline{\l}_r-1,\chiba_1,\chiba_s)$ 
is a $k$-basis of $\Rad k\groth(G)$.\finl

\bigskip

\subsection{Projective modules\label{subsection cartan}}
We now fix a conjugacy class $C$ of $p$-regular elements 
(i.e. $C \in G_{p'}/\ssim$). Let 
$$e_C=1_{\SC_{p'}(C)}=\sum_{D \in \SC_{p'}(C)/\!\sim} 1_D.$$
If necessary, $e_C$ will be denoted by $e_C^G$. 
If $H$ is a subgroup of $G$, then 
\equat\label{restriction ec}
\Res_H^G e_C^G = \sum_{D \in (C \cap H)/\!\sim_H} e_D^H.
\endequat

\begin{prop}\label{O}
Let $C \in G_{p'}/\ssim$. Then $e_C \in \OC_\pG\groth(G)$.
\end{prop}

\begin{proof}
Using Brauer's Theorem, we only need to prove that 
$\Res_N^G e_C^G \in \OC_\pG\groth(N)$ for every nilpotent subgroup $N$ 
of $G$. By \ref{restriction ec}, this amounts to prove the 
lemma whenever $G$ is nilpotent. So we assume that $G$ is nilpotent. 
Then $G=G_{p'} \times G_p$, and $G_p$ and $G_{p'}$ are subgroups of $G$. 
Moreover, $C \subset G_{p'}$ and $\SC_{p'}(G)=C \times G_p$. 
If we identify $K\groth(G)$ and 
$K\groth(G_{p'}) \otimes_K K\groth(G_p)$, 
we have $e_C^G = 1_C^{G_{p'}} \otimes_{\OC_\pG} e_1^{G_p}$. But, by \ref{formule ec}, 
we have that $e_C^{G_{p'}} \in \OC_\pG\groth(G_{p'})$. On the other hand, 
$e_1^{G_p} = 1_{G_p} \in \groth(G_p)$. The proof of the lemma is complete.
\end{proof}

\begin{coro}\label{primitif}
Let $C \in G_{p'}/\ssim$. Then $e_C$ is a primitive idempotent of $\OC_\pG\groth(G)$.
\end{coro}

\begin{proof}
By Proposition \ref{simples} (a), the number of primitive idempotents 
of $k\groth(G)$ is $|G_{p'}/\ssim|$. So the number of primitive 
idempotents of $\hat{\OC}_\pG\groth(G)$ is also $|G_{p'}/\ssim|$ 
(here, $\hat{\OC}_\pG$ denotes the completion of $\OC_\pG$ at its maximal 
ideal). 
Now, $(e_C)_{C \in G_{p'}/\!\sim}$ is a family of 
orthogonal idempotents of $\OC_\pG\groth(G)$ (see Proposition \ref{O}) and 
$1_G=\sum_{C \in G_{p'}/\!\sim} e_C$. 
The proof of the lemma is complete.
\end{proof}

Let $\eba_C \in k\groth(G)$ denote the reduction modulo 
$\pG\OC_\pG$ of $e_C$. Then it follows from \ref{factorisation} that 
\equat
d_\pG^k \eba_C = 1_{\SC_{p'}(C)}^k \in k\groth(kG) \simeq \Class_k^{p'}(G).
\endequat
Let $\PC_C=\OC_\pG\groth(G)e_C$ 
and $\PCba_C=k\groth(G)\eba_C$: they are indecomposable 
projective modules for $\OC_\pG\groth(G)$ and $k\groth(G)$ respectively. Then
$$\OC_\pG\groth(G)=\mathop{\oplus}_{C \in G_{p'}/\!\sim} \PC_C$$
$$k\groth(G)=\mathop{\oplus}_{C  \in G_{p'}/\!\sim} \bar{\PC}_C.\leqno{\text{and}}$$
Note also that 
\equat\label{dimension blocs}
\dim_k k\groth(G)\eba_C = {\mathrm{rank}}_{\OC_\pG} \OC_\pG\groth(G)e_C=
|\SC_{p'}(G)/\ssim|.
\endequat

\begin{prop}\label{cartan}
Let $C$ and $C'$ be two conjugacy classes of $p'$-regular elements 
of $G$. Then:
\begin{itemize}
\itemth{a} $[\bar{\PC}_{C}:\bar{\DC}_{C'}] = \begin{cases}
|\SC_{p'}(C)/\ssim| & \text{if }C=C',\\
0 & \text{otherwise.}
\end{cases}$.

\itemth{b} $\bar{\PC}_C/\Rad \bar{\PC}_C \simeq \bar{\DC}_C$.
\end{itemize}
\end{prop}

\begin{proof}
Let us first prove (a). By definition of $e_C$, we have
$$[K \otimes_{\OC_\pG} \PC_C]= \sum_{D \in \SC_{p'}(G)/\!\sim} [\DC_D].$$ 
Also, by definition of the decomposition map 
$\d_\pG : \groth(K\groth(G)) \to \groth(k\groth(G))$, we have 
$$\d_\pG [K \otimes_{\OC_\pG} \PC_C]=[\PCba_C].$$
So the result follows from these observations and from \ref{decomposition}. 
Now, (b) follows easily from (a). 
\end{proof}

\bigskip 

\subsection{More on the radical} 
Let $\Rad_p(G)$ denote the set of functions $f \in \OC_\pG\groth(G)$ 
whose restriction to $G_{p'}$ is zero. Note that $\Rad_p(G)$ is a 
direct summand of the $\OC_\pG$-module $\OC_\pG\groth(G)$. So, 
$k\Rad_p(G)=k \otimes_{\OC_\pG} \Rad_p(G)$ is a sub-$k$-vector space of 
$k\groth(G)$. 

\begin{prop}\label{description radical}
We have:
\begin{itemize}
\itemth{a} $\dim_k k\Rad_p(G) = |G/\ssim|-|G_{p'}/\ssim|$.

\itemth{b} $k\Rad_p(G)$ is the radical of $k\groth(G)$.
\end{itemize}
\end{prop}

\begin{proof}
(a) is clear. (b) follows from \ref{factorisation} and from Corollary 
\ref{simples}.
\end{proof}

\begin{coro}\label{exposant sylow}
Let $e$ be the number such that $p^e$ is the exponent of 
a Sylow $p$-subgroup of $G$. If $f \in \Rad k\groth(G)$, then $f^{p^e}=0$. 
\end{coro}

\begin{proof}
Let $e=e_p(G)$. If $f \in K\groth(G)$ and if $n \ge 1$, we denote by 
$f^{(n)} : G \to K$, $g \mapsto f(g^n)$. Then the map $K\groth(G) \to K\groth(G)$, 
$f \mapsto f^{(n)}$ is a morphism of $K$-algebras. Moreover (see for instance  
\cite[Corollary 12.10]{curtis}), we have
\equat\label{puissance stable}
\text{\it If $f \in \groth(G)$, then $f^{(n)} \in \groth(G)$.}
\endequat
Therefore, it induces a morphism of $k$-algebras $\theta_n : k\groth(G) \to k\groth(G)$. 
Now, let $F : k \groth(G) \to k\groth(G)$, $\l \otimes_\ZM f \mapsto \l^p \otimes_\ZM f$. 
Then $F$ is an injective endomorphism of the ring $k\groth(G)$. 
Moreover (see for instance \cite[Problem 4.7]{isaacs}), we have 
\equat
F \circ \th_p(f) = f^p
\endequat
for every $f \in k\groth(G)$. Since $F$ and $\th_p$ commute, we have 
$F^e \circ \th_{p^e}(f) = f^{p^e}$ for every $f \in k\groth(G)$. 
Therefore, if $\chi \in \Rad_p(G)$, we have 
$$\chiba^{p^e} = F^e(\overline{\chi^{(p^e)}}).$$
But, by hypothesis, $g^{p^e} \in G_{p'}$ for every $g \in G$. 
So, if $f \in \Rad_p(G)$, then $f^{(p^e)}=0$. Therefore, $\fba^{p^e}=0$. The corollary 
follows from this observation and from Proposition \ref{description radical}.
\end{proof}

\section{Principal block}

\medskip

If $C \in G_{p'}/\sim$, we denote by $\groth_\pG(G,C)$ the $\OC_\pG$-algebra 
$\OC_\pG\groth(G)e_C$. As an $\OC_\pG\groth(G)$-module, this is just $\PC_C$, 
but we want to study here its structure as a ring, so that is why we use 
a different notation. If $R$ is a commutative $\OC_\pG$-algebra, we set 
$R\groth_\pG(G,C)=R \otimes_{\OC_\pG} \groth_\pG(G,C)$. For instance, 
$k\groth_\pG(G,C)=k\groth(G)\eba_C$, and $K\groth_\pG(G,C)$ can be identified 
with the algebra of class functions on $\SC_{p'}(C)$. 

The algebra $\groth_\pG(G,1)$ (resp. $k\groth_\pG(G,1)$) will be called 
the {\it principal block} of $\OC_\pG\groth(G)$ (resp. $k\groth(G)$). 
The aim of this section is to construct an isomorphism 
$\groth_\pG(G,C) \simeq \groth_\pG(C_G(g),1)$, where $g$ is any element 
of $C$. We also emphasize the functorial properties of the principal 
block. 

\bigskip

\remark{tz ec} 
If $C \in G_{p'}/\ssim$ and if $z \in Z(G)$, then $t_z$ induces an isomorphism 
of algebras $\groth_\pG(G,C) \simeq \groth_\pG(G,z_{p'}^{-1} C)$ (see Remark 
\ref{tz unc}). Consequently, $\bar{t}_z$ induces an isomorphism of algebras 
$k\groth_\pG(G,C) \simeq k\groth_\pG(G,z^{-1} C)$.\finl

\bigskip

\subsection{Centralizers} 
Let $C \in G_{p'}/\ssim_G$. Let 
$\proj_C^G : K\groth(G) \to K\groth_\pG(G,C)$, 
$x \mapsto x e_C$ denote the canonical projection. We still denote by 
$\proj_C^G : \OC_\pG\groth(G) \to \groth_\pG(G,C)$, 
the restriction of $\proj_C^G$ and we denote by 
$\projbar_C^G : k\groth(G) \to k\groth_\pG(G,C)$ its reduction modulo 
$\pG\OC_\pG$. 

Let us now fix $g \in C$. 
It is well-known (and easy) that the map 
$C_G(g)_p/\ssim_{C_G(g)} \to \SC_{p'}(C)/\ssim_G$ that sends 
the $C_G(g)$-conjugacy class $D \in C_G(g)_p/\ssim_{C_G(g)}$ to the $G$-conjugacy 
class containing $gD$ is bijective. In particular,
\equat\label{cardinal}
|\SC_{p'}(C)/\ssim_G|=|C_G(g)_p/\ssim_{C_G(g)}|.
\endequat
Now, let $d_g^G : K\groth(G) \to K\groth(C_G(g))$ be the map defined by:
$$(d_g^G f)(h)=\begin{cases}
f(gh) & \text{if $h \in C_G(g)_p$,}\\
0 & \text{otherwise,}
\end{cases}$$
for all $f \in K\groth(G)$ and $h \in C_G(g)$. Then 
$d_g^G f \in K\groth_\pG(C_G(g),1)$. It must be noticed that 
\equat\label{composition}
d_g^G = \proj_1^{C_G(g)} \circ t_g^{C_G(g)} \circ \Res_{C_G(g)}^G 
= t_g^{C_G(g)} \circ \proj_g^{C_G(g)} \circ \Res_{C_G(g)}^G. 
\endequat
In particular, $d_g^G$ sends $\OC_\pG\groth(G)$ to $\groth_\pG(C_G(g),1)$. 
We denote by $\res_g : \groth_\pG(G,C) \to \groth_\pG(C_G(g),1)$ 
the restriction of $d_g^G$ to $\groth_\pG(G,C)$. Let 
$\ind_g : K\groth_\pG(C_G(g),1) \to K\groth_\pG(G,C)$ be the map 
defined by 
$$\ind_g f = \Ind_{C_G(g)}^G(t_{g^{-1}}^{C_G(g)} f)$$
for every $f \in K\groth_\pG(C_G(g),1)$. It is clear that 
$\ind_g f \in \groth_\pG(G,C)$ if $f \in \groth_\pG(C_G(g),1)$. 
Thus we have defined two maps
$$\res_g : \groth_\pG(G,C) \to \groth_\pG(C_G(g),1)$$
$$\ind_g : \groth_\pG(C_G(g),1) \to \groth_\pG(G,C).
\leqno{\text{and}}$$
We have:

\begin{theo}\label{isomorphismes}
If $g \in G_{p'}$, then $\res_g$ and $\ind_g$ are 
isomorphisms of $\OC_\pG$-algebras inverse to each other.
\end{theo}

\begin{proof}
We first want to prove that $\res_g \circ \ind_g$ is the identity. 
Let $f \in K\groth_\pG(C_G(g),1)$. Let $f' = t_{g^{-1}} f$ and let 
$x \in C_G(g)_p$. We just need to prove that 
$$(\Ind_{C_G(g)}^G f')(gx)=f'(gx).\leqno{(?)}$$
But, by definition, 
$$(\Ind_{C_G(g)}^G f')(gx)=\sum_{\substack{h \in [G/C_G(g)] \\ h(gx)h^{-1} \in C_G(g)}} 
f'(h(gx)h^{-1}).$$
Here, $[G/C_G(g)]$ denotes a set of representatives of $G/C_G(g)$. Since 
$f'$ has support in $gC_G(g)_p$, we have $f(h(gx)h^{-1}) \neq 0$ only if 
the $p'$-part of $h(gx)h^{-1}$ is equal to $g$, which happens if and only 
if $h \in C_G(g)$. This shows $(?)$. 

The fact that $\ind_g \circ \res_g$ is the identity can be proved similarly, or 
can be proved by using a trivial dimension argument. Since $\res_g$ is a morphism of 
algebras, we get that $\ind_g$ is also a morphism of algebras.
\end{proof}

\bigskip

\subsection{Subgroups of index prime to ${\boldsymbol{p}}$} 
If $H$ is a subgroup of $G$, then the restriction 
map $\Res_H^G$ sends $\groth_\pG(G,1)$ to $\groth_\pG(H,1)$ 
(indeed, by \ref{restriction ec}, we have $\Res_H^G e_1^G=e_1^H$).

\begin{theo}\label{facteur direct}
If $H$ is a subgroup of $G$ of index prime to $p$, then 
$\Res_H^G : \groth_\pG(G,1) \to \groth_\pG(H,1)$ 
is a split injection of $\OC_\pG$-modules.
\end{theo}

\begin{proof}
Let us first prove that $\Res_H^G$ is injective. For this, we only need 
to prove that the map $\Res_H^G : K\groth_\pG(G,1) \to K\groth_\pG(H,1)$. 
But $K\groth_\pG(G,1)$ is the space of functions whose support is contained 
in $G_p$. Since the index of $H$ is prime to $p$, every conjugacy class 
of $p$-elements of $G$ meets $H$. This shows that $\Res_H^G$ is injective. 

In order to prove that it is a split injection, we only need to prove that 
the $\OC_\pG$-module 
$\groth_\pG(H,1)/\Res_H^G(\groth_\pG(G,1))$ is torsion-free. 
Let $\pi$ be a generator of the ideal $\pG\OC_\pG$. Let $\g \in \groth_\pG(G,1)$ 
and $\eta \in \groth_\pG(H,1)$ be such that $\pi \eta = \Res_H^G \g$. 
We only need to prove that $\g/\pi \in \groth_\pG(G,1)$. 
By Brauer's Theorem, it is sufficient to show that, for any nilpotent 
subgroup $N$ of $G$, we have $\Res_N^G \g \in \pi\OC_\pG\groth(N)$. 

So let $N$ be a nilpotent subgroup. We have $N=N_p \times N_{p'}$ and, since 
the index of $H$ in $G$ is prime to $p$, we may assume that $N_p \subset H$. 
Since $\Res_N^G \psi \in \groth_\pG(N,1) = 
\OC_\pG\groth(N_p) \otimes_{\OC_\pG} e_1^{N_{p'}}$, we have  
\eqna
\Res_N^G \g &=& (\Res_{N_p}^G \g) \otimes_{\OC_\pG} e_1^{N_{p'}} \\
&=& (\pi \Res_{N_p}^H \eta) \otimes_{\OC_\pG} e_1^{N_{p'}}  \in \pi\OC_\pG \groth(N),
\endeqna
as expected.
\end{proof}

\begin{coro}\label{injection}
If $H$ is a subgroup of $G$ of index prime to $p$, then the map 
$\Resbar_H^G : k\groth_\pG(G,1) \to k\groth_\pG(H,1)$ 
is an injective morphism of $k$-algebras.
\end{coro}

\begin{coro}\label{controle de fusion}
If $H$ is a subgroup of $G$ of index prime to $p$ which controls 
the fusion of $p$-elements, then 
$\Res_H^G : \groth_\pG(G,1) \to \groth_\pG(H,1)$ 
is an isomorphism of $\OC_\pG$-algebras.
\end{coro}

\begin{proof}
In this case, $\dim_K K\groth_\pG(G,1)=\dim_K K\groth_\pG(H,1)$, so 
the result follows from Corollary \ref{injection}.
\end{proof}

\bigskip

\example{sylow abelien}
Let $P$ be a Sylow $p$-subgroup of $G$ and assume in this example that 
$P$ is abelian. Then $N_G(P)$ controls the fusion of $p$-elements. It then follows 
from Corollary \ref{controle de fusion} that the restriction from $G$ to $N_G(P)$ 
induces isomorphisms of algebras 
$\groth_\pG(G,1) \simeq \groth_\pG(N_G(P),1)$ and 
$k\groth_\pG(G,1) \simeq k\groth_\pG(N_G(P),1)$. In particular, 
$\loewy_p(G,1)=\loewy_p(N_G(P),1)$.\finl

\bigskip

\example{semi direct} 
Let $N$ be a $p'$-group, let $H$ be a group acting on $N$ and let 
$G=H \ltimes N$. Then $H$ is of index prime to $p$ and 
controls the fusion of $p$-elements of $G$. So $\Res_H^G$ induces 
isomorphisms of algebras 
$\groth_\pG(G,1) \simeq \groth_\pG(H,1)$ and 
$k\groth_\pG(G,1) \simeq k\groth_\pG(H,1)$. In particular, 
$\loewy_p(G,1)=\loewy_p(H,1)$.\finl

\bigskip

\subsection{Quotient by a normal ${\boldsymbol{p'}}$-subgroup} 
Let $N$ be a normal subgroup of $G$. Let $\pi : G \to G/N$ 
denote the canonical morphism. Then the morphism of 
algebras $\Res_\pi : \groth_\pG(G/N) \to \groth_\pG(G)$ induces a 
morphism of algebras $\Res_\pi^{(1)} : \groth_\pG(G/N,1) \to \groth_\pG(G,1)$, 
$f \mapsto (\Res_\pi f)e_1^G$. Note that $\Res_\pi^{(1)} e_1^{G/N} = e_1^G$. 
We denote by $\Resbar_\pi^{(1)} : k\groth_\pG(G/N,1) \to k\groth_\pG(G,1)$ the 
morphism induced by $\Res_\pi^{(1)}$. 
Then:

\begin{theo}\label{p' quotient}
With the above notation, we have:
\begin{itemize}
\itemth{a} $\Res_\pi^{(1)}$ is a split injection of $\OC_\pG$-modules.

\itemth{b} If $N$ is prime to $p$, then $\Res_\pi^{(1)}$ is an isomorphism.
\end{itemize}
\end{theo}

\begin{proof}
(a) The injectivity of $\Res_\pi^{(1)}$ follows from the fact that $(G/N)_p=G_pN/N$. 
Now, let $I$ denote the image of $\Res_\pi^{(1)}$. Since $\Res_\pi(\OC_\pG\groth(G/N))$ 
is a direct summand of $\OC_\pG\groth(G)$, we get that 
$\Res_\pi(\groth_\pG(G/N,1))=(\Res_\pi^{(1)} e_1^{G/N}) \Res_\pi(\OC_\pG\groth(G/N))$ 
is a direct summand of $\OC_\pG\groth(G)$. Since $I=e_1^G\Res_\pi(\groth_\pG(G/N,1))$ 
and $e_1^G=e_1^G \Res_\pi(e_1^{G/N})$, 
we get that $I=e_1^G \Res_\pi(\OC_\pG\groth(G/N))$ is a direct summand 
of $\OC_\pG\groth(G)$, as desired.

\medskip

(b) now follows from (a) and from the fact that the map $\pi$ induces 
a bijection between $G_p/\ssim_G$ and $(G/N)_g/\ssim_{G/N}$ whenever 
$N$ is a normal $p'$-subgroup.
\end{proof}

\bigskip

\section{Some invariants}

\medskip

We introduce in this section some numerical invariants of the $k$-algebra 
$k\groth(G)$ (more precisely, of the algebras $k\RC_\pG(G,C)$): Loewy length, 
dimension of the $\Ext$-groups. 

\bigskip

\subsection{Loewy length} 
If $C \in G_{p'}/\ssim$, we denote by $\loewy_p(G,C)$ the Loewy length of 
the $k$-algebra $k\groth_\pG(G,C)$. Then, by definition, we have 
\equat\label{debile}
\loewy_p(G)=\max_{C \in G_{p'}/\!\sim} \loewy_p(G,C).
\endequat
On the other hand, by Theorem \ref{isomorphismes}, we have 
\equat\label{egalite loewy}
\text{\it If $C \in G_{p'}/\ssim$ and if $g \in C$, then 
$\loewy_p(G,C)=\loewy_p(C_G(g),1)$.}
\endequat
The following bound on the Loewy length of $k\groth(G)$ is obtained 
immediately from \ref{dimension blocs} and \ref{cardinal}:
\equat\label{borne}
\loewy_p(G) \le \max_{C \in G_{p'}/\!\sim} |\SC_{p'}(C)/\ssim|=
\max_{g \in G_{p'}} |C_G(g)_p /\ssim_{C_G(g)}|.
\endequat
We set $S_p(G)=\DS{\max_{C \in G_{p'}/\!\sim}} |\SC_{p'}(C)/\ssim|$. 

\bigskip

\example{borne mauvaise}
The inequality \ref{borne} might be strict. Indeed, if 
$G=\ZM/2\ZM \times \ZM/2\ZM$, then 
$\loewy_2(G)=3 < 4=S_2(G)$.\finl

\bigskip

\example{borne bonne}
If $S_p(G)=2$, then $\loewy_p(G)=2$. Indeed, in this case, we have that $p$ divides 
$|G|$, so $k\groth(G)$ is not semisimple by Corollary \ref{critere semisimple}, 
so $\loewy_p(G) \ge 2$. The result then follows from \ref{borne}.\finl

\bigskip

\subsection{${\mathbf{Ext}}$-groups} 
If $i \ge 0$ and if $C \in G_{p'}/\ssim$, we set 
$$\ext_p^i(G,C)=\dim_\gfp \Ext_{k\groth(G)}^i(\DCba_C,\DCba_C).$$
Note that $\ext_p^i(G,C)=\dim_\gfp \Ext_{k\groth_\pG(G,C)}^i(\DCba_C,\DCba_C)$. 
So, if $g \in C$, it follows from Theorem \ref{isomorphismes} that 
\equat\label{ext}
\ext_p^i(G,C)=\ext_p^i(C_G(g),1).
\endequat

\bigskip

\subsection{Subgroups, quotients} 
The next results follows respectively from Corollaries \ref{injection}, 
\ref{controle de fusion} and from Theorem \ref{p' quotient}:

\begin{prop}\label{borne meilleure}
Let $H$ be a subgroup of $G$ of index prime to $p$ and let $N$ be a normal 
subgroup of $G$. Then:
\begin{itemize}
\itemth{a} $\loewy_p(G,1) \leqslant \loewy_p(H,1)$.

\itemth{b} If $H$ controls the fusion of $p$-elements, then 
$\loewy_p(G,1) = \loewy_p(H,1)$ and $\ext_p^i(G,1)=\ext_p^i(H,1)$ for every $i \ge 0$.

\itemth{c} $\loewy_p(G/N,1) \leqslant \loewy_p(G,1)$.

\itemth{d} If $|N|$ is prime to $p$, then 
$\loewy_p(G,1) = \loewy_p(H,1)$ and $\ext_p^i(G,1)=\ext_p^i(H,1)$ for every $i \ge 0$.
\end{itemize}
\end{prop}

%
%
\bigskip

\subsection{Direct products} 
We study here the behaviour of the invariants $\ell_p(G,C)$ and $\ext_p^1(G,C)$ 
with respect to taking direct products. We first recall the following 
result on finite dimensional algebras:

\begin{prop}\label{tenseur}
Let $A$ and $B$ be two finite dimensional $k$-algebras. Then:
\begin{itemize}
\itemth{a} $\Rad(A \otimes_k B) = A \otimes_k (\Rad B) + (\Rad A) \otimes_k B$. 

\itemth{b} If $A/\Rad A \simeq k$ and $B/\Rad B \simeq k$, then 
$$\Rad(A \otimes_k B)/\Rad(A \otimes_k B)^2 \simeq (\Rad A)/(\Rad A)^2 \oplus 
(\Rad B)/(\Rad B)^2.$$
\end{itemize}
\end{prop}

\begin{proof}
(a) is proved for instance in \cite[Proof of 10.39]{curtis}. 
Let us now prove (b). 
Let $\th : (\Rad A) \oplus (\Rad B) \to 
\Rad(A \otimes_k B)/\Rad(A \otimes_k B)^2$, 
$a \oplus b \mapsto \overline{a \otimes_k 1 + 1 \otimes_k b}$. 
By (a), $\th$ is surjective and 
$(\Rad A)^2 \oplus (\Rad B)^2$ is contained in the kernel of $\th$. 
Now the result follows from dimension reasons (using (a)).
\end{proof}

\begin{prop}\label{direct product}
Let $G$ and $H$ be two finite groups and let $C \in G_{p'}/\ssim$ and 
$D \in H_{p'}/\ssim$. Then 
$$\loewy_p(G \times H, C \times D) =\loewy_p(G,C) + \loewy_p(H,D) - 1$$
$$\ext_p^1(G\times H,C \times D)=\ext_p^1(G,C) + \ext_p^1(H,D).\leqno{\text{and}}$$
\end{prop}

\begin{proof}
Write $A=k\groth_\pG(G,C)$ and $B=k\groth_\pG(H,D)$. 
It is easily checked that $k\groth_\pG(G \times H,C \times D)= A \otimes_k B$. 
So the first equality follows from Propositon \ref{tenseur} (a) and 
from the commutativity of $A$ and $B$. 
Moreover $A/(\Rad A)\simeq k$ and $B/(\Rad B) \simeq k$. 
In particular $\dim_k \Ext_A^1(A/\Rad A,A/\Rad A)=\dim_k (\Rad A)/(\Rad A)^2$. 
So the second equality follows from Proposition \ref{tenseur} (b).
\end{proof}

\bigskip

\subsection{Abelian groups} 
We compute here the invariants $\ell_p(G,1)$ and $\ext_p^1(G,1)$ whenever $G$ 
is abelian. If $G$ is abelian, then there is a (non-canonical) isomorphism 
of algebras $k\groth(G) \simeq kG$. 

Let us first start with the cyclic case:
\equat\label{cyclique}
\text{\it if $G$ is cyclic, then $\loewy_p(G)=|G|_p+1$ and 
$\ext_p^1(G,1)=\begin{cases} 1 & \text{\it if $p$ divides $|G|$,} \\ 
0 & \text{\it otherwise.}\end{cases}$}
\endequat
Therefore, by Proposition \ref{direct product}, we have: 
if $G_1$, \dots, $G_n$ are cyclic, then 
\equat\label{abelien loewy}
\loewy_p(G_1 \times \dots \times G_n) = |G_1|_p + \dots + |G_n|_p -n+1.
\endequat
and
\equat\label{abelien ext}
\ext_p^1(G_1 \times \dots \times G_n) = |\{1 \le i \le n~|~\text{\it $p$ divides $G_i$}\}|.
\endequat

\medskip

\section{The symmetric group}

\medskip

\def\grothbar{\overline{\groth}}

In this section, and only in this section, we fix a non-zero natural number 
$n$ and a prime number $p$ and we assume that $G=\SG_n$, that $\OC=\ZM$ and 
that $\pG=p\ZM$. Let $\gfp=k$. It is well-known that $\QM$ and $\gfp$ are splitting 
fields for $\SG_n$. For simplification, we set $\RC_n=\RC(\SG_n)$ and 
$\overline{\groth}_n=\gfp\groth(\SG_n)$. We investigate further the structure 
of $\grothbar_n$. This is a continuation 
of the work started in \cite{bonnafe sym} in which 
the description of the descending Loewy series of 
$\grothbar_n$ was obtained.

We first introduce some notation. Let $\Part(n)$ denote the set of partitions 
of $n$. If $\l=(\l_1,\dots,\l_r) \in \Part(n)$ and if $1 \le i \le n$, we denote 
by $r_i(\l)$ the number of occurences of $i$ as a part of $\l$. We set
$$\pi_p(\l)=\sum_{i = 1}^n \Bigl[\frac{r_i(\l)}{p}\Bigr]$$
where, for $x \in \RM$, $x \ge 0$, we denote by $[x]$ the unique natural number 
$m \ge 0$ such that $m \le x < m+1$. Note that $\pi_p(\l) \in \{0,1,2,\dots,[n/p]\}$ 
and recall that $\l$ is {\it $p$-regular} (resp. {\it $p$-singular}) if and only if 
$\pi_p(\l)=0$ (resp. $\pi_p(\l) \ge 1$). 
We denote by $\SG_\l$ the Young subgroup canonically isomorphic to 
$\SG_{\l_1} \times \dots \times \SG_{\l_r}$, by $1_\l$ the trivial character of $\SG_\l$, 
and by $c_\l$ an element of $\SG_\l$ with only $r$ orbit in $\{1,2,\dots,n\}$. 
Let $C_\l$ denote the conjugacy class of $c_\l$ in $\SG_n$. Then the map 
$\Part(n) \to \SG_n/\ssim$, $\l \mapsto C_\l$ is a bijection. Let 
$W(\l)=N_{\SG_n}(\SG_\l)/\SG_\l$. Then 
\equat\label{norm sn}
W(\l) \simeq \prod_{i = 1}^n \SG_{r_i(\l)}.
\endequat
In particular, $\pi_p(\l)$ is the {\it $p$-rank} of $W(\l)$, where the $p$-rank of a finite 
group is the maximal rank of an elementary abelian subgroup. 
Now, we set $\ph_\l=\Ind_{\SG_\l}^{\SG_n} 1_\l$. An old result of Frobenius says 
that
\equat\label{base}
\text{\it $(\ph_\l)_{\l \in \Part(n)}$ is a $\ZM$-basis of $\groth_n$}
\endequat
(see for instance \cite[Theorem 5.4.5 (b)]{geck}). Now, if $i \ge 1$, let 
$$\Part_p^{\geqslant i}(n)=\{\l \in \Part(n)~|~\pi_p(\l) \ge i\}$$
$$\Part_p^i(n)=\{\l \in \Part(n)~|~\pi_p(\l) = i\}.\leqno{\mathrm{and}}$$
Then, by \cite[Theorem A]{bonnafe sym}, we have
\equat\label{rad sym}
(\Rad \grothbar_n)^i = \mathop{\oplus}_{\l \in \Part_p^{(i)}(n)} \gfp \phba_\l.
\endequat

Let $\Part_{p'}(n)$ denote the set of partitions of $n$ 
whose parts are prime to $p$. Then the map $\Part_{p'}(n) \to G_{p'}/\ssim$, 
$\l \mapsto C_\l$ is bijective. We denote by $\t_{p'}(\l)$ the unique partition of $n$ 
such that $(c_\l)_{p'} \in C_{\t_{p'}(\l)}$. If $\l=(\l_1,\l_2,\dots,\l_r)$, 
the partition $\t_{p'}(\l)$ is obtained as follows. Let 
$$\l'=(\underbrace{(\l_1)_{p'},\dots,(\l_1)_{p'}}_{(\l_1)_p~\mathrm{times}},\dots,
\underbrace{(\l_r)_{p'},\dots,(\l_r)_{p'}}_{(\l_r)_p~\mathrm{times}}).$$
Then $\t_{p'}(\l)$ is obtained from $\l'$ by reordering the parts. 
The map $\t_{p'} : \Part(n) \to \Part_{p'}(n)$ is obviously surjective. 
If $\l \in \Part_{p'}(n)$, we set for simplification 
$\groth_{n,p}(\l)=\groth_{p\ZM}(\SG_n,C_\l)$ and 
$\grothbar_n(\l)=\gfp\groth_{p\ZM}(\SG_n,C_\l)$. In other words, 
$$\ZM_{p\ZM}\groth_n = \mathop{\oplus}_{\l \in \Part_{p'}(n)} \groth_{n,p}(\l)$$
$$\grothbar_n=\mathop{\oplus}_{\l \in \Part_{p'}(n)} \grothbar_n(\l)\leqno{\text{and}}$$
are the decomposition of $\ZM_{p\ZM}\groth_n$ and $\grothbar_n$ as a sum of blocks. 
We now make the result \ref{rad sym} more precise:

\begin{prop}\label{rad sym precis}
If $\l \in \Part_{p'}(n)$ and if $i \ge 0$, then 
$$\dim_\gfp \bigl(\Rad \grothbar_n(\l)\bigr)^i=|\t_{p'}^{-1}(\l) \cap 
\Part_p^{\geqslant i}(n)|.$$
\end{prop}

\def\m{\mu}

\begin{proof}
If $\l$ and $\m$ are two partitions of $n$, we write $\l \subset \m$ if $\SG_\l$ 
is conjugate to a subgroup of $\SG_\m$. This defines an order on $\Part(n)$. 
On the other hand, if $d \in \SG_n$, we denote by $\l \cap \lexp{d}{\m}$ 
the unique partition $\nu$ of $n$ such that $\SG_\l \cap \lexp{d}{\SG_\m}$ 
is conjugate to $\SG_\nu$. Then, by the Mackey formula for tensor product 
(see for instance \cite[Theorem 10.18]{curtis}), we have 
$$\ph_\l\ph_\m = \sum_{d \in [\SG_\l\backslash \SG_n / \SG_\m]} \ph_{\l \cap \lexp{d}{\m}}.
\leqno{(1)}$$
Here, $[\SG_\l\backslash \SG_n / \SG_\m]$ denotes a set of representatives 
of the $(\SG_\l,\SG_\m)$-double cosets in $\SG_n$. This shows that, if we fixe 
$\l_0 \in \Part(n)$, then $\oplus_{\l \subset \l_0} \ZM\ph_\l$ and 
$\oplus_{\l \subsetneq \l_0} \ZM\ph_\l$ are sub-$\groth(G)$-module of $\groth(G)$. 
We denote by $\DC_\l^\ZM$ the quotient of these two modules. Then 
$$K \otimes_\ZM \DC_\l^\ZM \simeq \DC_{C_\l}.\leqno{(2)}$$
This follows for instance from \cite[Proposition 2.4.4]{geck}. Consequently, 
$$k \otimes_\ZM \DC_\l^\ZM \simeq \DCba_{C_\l}.\leqno{(3)}$$
It then follows from Proposition \ref{isomorphisme} that 
$$\text{\it $k \otimes_\ZM \DC_\l^\ZM\simeq k \otimes_\ZM \DC_\m^\ZM~$ 
if and only if $~\t_{p'}(\l)=\t_{p'}(\m)$.}\leqno{(4)}$$
Now the Theorem follows from easily from (3), (4) and \ref{rad sym}.
\end{proof}

Now, if $\l \in \Part_{p'}(n)$, then $C_{\SG_n}(w_\l)$ contains a 
normal $p'$-subgroup $N_\l$ such that $C_{\SG_n}(w_\l)/N_\l \simeq W(\l)$. 
We denote by $1^n$ the partition $(1,1,\dots,1)$ of $n$. 
It follows from Theorem \ref{isomorphismes} and Theorem \ref{p' quotient} 
that 
\equat\label{iso rn}
\groth_{n,p}(\l) \simeq \groth_{p\ZM}(W(\l),1) \simeq \mathop{\otimes}_{i=1}^n 
\groth_{r_i(\l),p}(1^{r_i(\l)})
\endequat
and
\equat\label{iso rnbar}
\grothbar_n(\l) \simeq \grothbar(W(\l),1) \simeq \mathop{\otimes}_{i=1}^n 
\grothbar_{r_i(\l)}(1^{r_i(\l)}).
\endequat
We denote by $\Log_p n$ the real number $x$ such that $p^x=n$. Then:

\begin{coro}\label{ext 1}
If $\l \in \Part_{p'}(n)$, then 
$$\ext_p^1(\SG_n,C_\l)=\sum_{i=1}^n [\Log_p r_i(\l)]$$
$$\ell_p(\SG_n,C_\l)=\pi_p(\l)+1.\leqno{\text{and}}$$
\end{coro}

\begin{proof}
By \ref{iso rnbar} and by Proposition \ref{direct product}, both equalities 
need only to be proved whenever $\l=(1^n)$. So we assume that $\l=(1^n)$. 

\smallskip

Let us show the first equality. By Proposition \ref{rad sym precis}, 
we are reduced to show that $|\t_{p'}^{-1}(1^n) \cap \Part_p^1(n)|=[\Log_p n]$. 
Let $r=[\Log_p n]$. In other words, we 
have $p^r \le n < p^{r+1}$. If $1 \le i \le r$, write 
$n-p^i=\sum_{j=0}^r a_{ij} p^j$ with $0 \le a_{ij} < p-1$ (the $a_{ij}$'s are 
uniquely determined). Let 
$$\l(i)=(\underbrace{p^r,\dots,p^r}_{a_{ir}~\mathrm{times}},\dots,
\underbrace{p^i,\dots,p^i}_{a_{ir}~\mathrm{times}}
\underbrace{p^{i-1},\dots,p^{i-1}}_{(p+a_{i-1,r})~\mathrm{times}},
\underbrace{p^{i-2},\dots,p^{i-2}}_{a_{i-2,r}~\mathrm{times}},\dots,
\underbrace{1,\dots,1}_{a_{0r}~\mathrm{times}}).$$
The result will follow from the following equality
$$\t_{p'}^{-1}(1^n) \cap \Part_p^1(n)=\{\l(1),\l(2),\dots,\l(r)\}.\leqno{(*)}$$
So let us now prove $(*)$. Let $I=\{\l(1),\l(2),\dots,\l(r)\}$. It is clear 
that $I \subset \t_{p'}^{-1}(1^n) \cap \Part_p^1(n)$. Now, let 
$\l \in \t_{p'}^{-1}(1^n) \cap \Part_p^1(n)$. Then there exists a unique 
$i \in \{1,2,\dots,r\}$ such that $r_{p^{i-1}}(\l)\ge p$. Moreover, 
$r_{p^{i-1}}(\l) < 2p$. So, if we set $r_{p^j}'=r_{r_j}(\l)$ if $j \neq i-1$ 
and $r_{p^{i-1}}'=r_{p^{i-1}}(\l)-p$, we get that 
$0 \le r_{p^j}' \le p-1$ and $n-p^i=\sum_{j=0}^r r_{p^j}' p^j$. 
This shows that $r_{p^j}' = a_{ij}$, so $\l=\l(i)$.

\medskip

Let us now show the second equality fo the Corollary. By Proposition \ref{rad sym precis}, 
we only need to show that 
$|\t_{p'}^{-1}(1^n) \cap \Part_p^{[n/p]}(n)|\ge 1$. But in fact, it is clear that 
$\t_{p'}^{-1}(1^n) \cap \Part_p^{[n/p]}(n)=\{1^n\}$. 
\end{proof}

\begin{coro}\label{ext 1 principal}
We have 
$$\dim_\gfp \bigl(\Rad \grothbar_n(1^n)\bigr)^{[n/p]}=1$$
$$\dim_\gfp \Ext_{\grothbar_n}^1(\DCba_{1^n},\DCba_{1^n})=[\Log_p n].\leqno{\text{and}}$$
In particular, $\ell_p(\SG_n,1)=\ell_p(\SG_n)=[n/p]$. 
\end{coro}

\begin{proof}
This is just a particular case of the previous corollary. The first equality 
has been obtained in the course of the proof of the previous corollary.
\end{proof}

\bigskip

\section{Dihedral groups}

\medskip

Let $n \ge 1$ and $m \ge 0$ be two natural numbers. We assume 
in this section, and only in this subsection, that $G=D_{2^n(2m+1)}$ is 
the dihedral group of order $2^n(2m+1)$ and that $p=2$. 

\bigskip

\begin{prop}
If $n \ge 1$ and $m \ge 0$ are natural numbers, then 
$$\loewy_2(D_{2^n(2m+1)},1)=\begin{cases}
2 & \text{if }n=1,\\
3 & \text{if }n=2,\\
2^{n-2}+1 & \text{if }n \ge 3.\\
\end{cases}$$
and
$$\ext_2^1(D_{2^n(2m+1)},1)=\begin{cases}
1 & \text{if }n=1,\\
2 & \text{if }n=2,\\
3 & \text{if }n \ge 3.
\end{cases}$$
\end{prop}

\begin{proof}
Let $N$ be the normal subgroup of $G$ of order $2m+1$. Then 
$G \simeq D_{2^n} \ltimes N$. So, by Proposition \ref{borne meilleure} (d), we may, and we 
will, assume that $m=0$. If $n=1$ or $2$ the the result is easily checked. Therefore, 
we may, and we will, assume that $n \ge 3$.

Write $h=2^{n-1}$. We have 
$$G=<s,t~|~s^2=t^2=(st)^h=1>.$$
Let $H=<st>$ and $S=<s>$. Then $|H|=2^{n-1}=h$ 
and $G=S \ltimes H$. We fix a primitive $h$-th root of unity $\z \in \OC^\times$. 
If $i \in \ZM$, we denote by $\xi_i$ the unique linear character 
of $H$ such that $\xi_i(st)=\z^i$. Then $\Irr H = \{\xi_0,\xi_1,\dots,\xi_{h-1}\}$, 
and $\xi_0=1_H$.

Since $n \ge 3$, $h$ is even and, if we write $h=2h'$, then $h'=2^{n-2}$ 
is also even. For $i \in \ZM$, we set 
$$\chi_i=\Ind_H^G \xi_i.$$
It is readily seen that $\chi_i=\chi_{-i}$, that $\chi_{i+h}=\chi_i$ and that 
\equat\label{produit diedral}
\chi_i \chi_j = \chi_{i+j} + \chi_{i-j}.
\endequat
Let $\e$ (resp. $\e_s$, resp. $\e_t$) be the unique linear 
character of order $2$ such that $\e(st)=1$ (resp. $\e_s(s)=1$, resp. 
$\e_t(t)=1$). Then 
$$\chi_0 = 1_G + \e,$$
$$\chi_{h'}=\e_s+\e_t,$$
and, if $h'$ does not divide $i$, 
$$\chi_i \in \Irr G.$$
Moreover, $|\Irr G|=h'+3$ and
$$\Irr G = \{1_G,\e,\e_s,\e_t,\chi_1,\chi_2,\dots,\chi_{h'-1}\}.$$
Finally, note that 
\equat\label{mult epsilon}
\e_s \chi_i = \e_t\chi_i =\chi_{i+h'}.
\endequat

Let us start by finding a lower bound for $\loewy_2(G)$. 
First, notice that the following equality holds: for all $i$, $j \in \ZM$ 
and every $r \ge 0$, we have 
\equat\label{puissance facile}
(\chiba_i+\chiba_j)^{2^r}=\chiba_{2^r i}+\chiba_{2^r j}.
\endequat
\begin{quotation}
\begin{proof}[Proof of \ref{puissance facile}]
Recall that $\chiba_i$ denotes the image of $\chi_i$ in 
$k\groth(G)$. We proceed by induction on $r$. The case $r=0$ is trivial. 
The induction step is an immediate consequence of \ref{produit diedral}.
\end{proof}
\end{quotation}

Note also the following fact (which follows from Example \ref{p group}):
\equat\label{radical exemples}
\text{\it If $i \in \ZM$, then $\chiba_i \in \Rad k\groth(G)$.}
\endequat

Therefore,
\equat\label{borne inf}
\loewy_2(G) \ge 2^{n-2}+1.
\endequat
\begin{quotation}
\begin{proof}[Proof of \ref{borne inf}]
By \ref{puissance facile}, we have immediately that 
$(\chiba_0+\chiba_1)^{2^{n-2}}=\chiba_0 + \chiba_{h'} \neq 0$ 
and, by \ref{radical exemples}, 
$\chiba_0+\chiba_1 \in \Rad k\groth(G)$.
\end{proof}
\end{quotation}
\bigskip

By Example \ref{p group}, we have 
\equat\label{base radical}
\text{\it $(\bar{1}_G+\epsba_s,\chiba_0,\chiba_1,\dots,\chiba_{h'})$ 
is a $k$-basis of $\Rad k\groth(G)$.}
\endequat
By \ref{mult epsilon} and \ref{produit diedral}, we get that 
\equat\label{base 2}
\text{\it $(\chiba_i+\chiba_{i+2})_{0 \le i \le h'-2}$ is a $k$-basis of 
$(\Rad k\groth(G))^2$.}
\endequat
This shows that $\ext_p^1(G)=3$, as expected. 
It follows that, if $n \ge 3$ and $2 \le i \le 2^{n-2}+1$, then 
\equat\label{diedral 2}
\dim_k\bigl(\Rad k\groth(D_{2^n})\bigr)^i = 2^{n-2}+1-i
\endequat
\begin{quotation}
\begin{proof}[Proof of \ref{diedral 2}]
Let $d_i = \dim_k\bigl(\Rad k\groth(D_{2^n})\bigr)^i$. By \ref{base 2}, 
we have $d_2 = 2^{n-2}-1$. By \ref{borne inf}, we have $d_{2^{n-2}} \ge 1$. 
Moreover, $d_1 > d_2 > d_3 > \dots$ So the proof of \ref{diedral 2} 
is complete.
\end{proof}
\end{quotation}
In particular, we get:
\equat\label{dernier}
\text{\it If $n \ge 3$, then 
$\bigl(\Rad k\groth(D_{2^n})\bigr)^{2^{n-2}}=
k(\bar{1}_{D_{2^n}} + \epsba + \epsba_s + \epsba_t)$.}
\endequat
and $\loewy_2(D_{2^n})=2^{n-2}+1$, as expected.
\end{proof}

\bigskip

\section{Examples}

\bigskip

For $0 \le i \le \loewy_p(G)-1$, we set $d_i=\dim_k (\Rad k\groth(G))^i$. 
Note that $d_0=|G/\ssim|$ and $d_0-d_1=|G_{p'}/\ssim|$. 
In this section, 
we give tables containing the values $\loewy_p(G)$, $\loewy_p(G,1)$, $S_p(G)$, 
$\ext_p^1(G,1)$ and the sequence ($d_0$, $d_1$, $d_2$,\dots) 
for various groups. These computations have been made using 
{\tt GAP3} \cite{GAP}.

These computations show that, if $G$ satisfies at least 
one of the following conditions:
\begin{quotation}
\noindent (1) $|G| \le 200$;

\noindent (2) $G$ is a subgroup of $\SG_8$;

\noindent (3) $G$ is one of the groups contained in the next tables;
\end{quotation}
then $\loewy_p(G,1)=\loewy_p(N_G(P),1)$ 
(here, $P$ denotes a Sylow $p$-subgroup of $G$). Note also that this 
equality holds if $P$ is abelian (see Example \ref{sylow abelien}). 

\medskip

\begin{quotation}
\noindent{\bf Question.} {\it Is it true that $\loewy_p(G,1)=\loewy_p(N_G(P),1)$?}
\end{quotation}

\medskip

The first table contains the datas for the the exceptional 
Weyl groups, the second table is for the alternating groups $\AG_n$ 
for $5 \le n \le 12$, the third table is for some small finite 
simple groups, and the last table is for the groups $P\!S\!L(2,q)$ 
for $q$ a prime power $\le 27$. 

$$\begin{array}{|c|c
|@{\!}c@{\!}||@{\!\hspace{0.15cm}}c@{\!\hspace{0.15cm}}|@{\!\hspace{0.15cm}}c@{\!\hspace{0.15cm}}|l|c|c|}
\hline
\vphantom{\DS{\frac{A}{B}}} 
G & |G| & \quad p\quad & \loewy_p(G) & S_p(G) & d_0, d_1, d_2, \dots & 
\loewy_p(G,1) & \ext_p^1(G,1) \\
\hline
\hline
\vertical W(E_6) & 51840     & 2 & 5 & 10 & 25, 19, 9, 3, 1 & 5 & 3 \\
\vertical        & 2^7.3^4.5 & 3 & 4 & 5  & 25, 13, 4, 1    & 4 & 2 \\
\vertical        &           & 5 & 2 & 2  & 25, 2           & 2 & 1 \\
\hline
\vertical W(E_7) & 2903040        & 2 & 7 & 24& 60, 52, 35, 18, 7, 3, 1 & 7 & 4 \\
\vertical        & 2^{10}.3^4.5.7 & 3 & 4 & 5 & 60, 30, 8, 2            & 4 & 2 \\
\vertical        &                & 5 & 2 & 2 & 60, 6                   & 2 & 1 \\
\vertical        &                & 7 & 2 & 2 & 60, 2                   & 2 & 1 \\
\hline
\vertical W(E_8) & 696729600        & 2 & 8 & 32& 112, 100, 68, 36, 17, 7, 3, 1 & 8 & 5 \\
\vertical        & 2^{14}.3^5.5^2.7 & 3 & 5 & 8 & 112, 65, 24, 7, 2              & 5 & 2 \\
\vertical        &                  & 5 & 3 & 3 & 112, 17, 2                     & 3 & 1 \\
\vertical        &                  & 7 & 2 & 2 & 112, 4                         & 2 & 1 \\
\hline
\vertical W(F_4) & 1152    & 2 & 5 & 14& 25, 21, 12, 4, 1 & 5 & 4 \\
\vertical        & 2^7.3^2 & 3 & 3 & 4 & 25, 11, 2        & 3 & 2 \\
\hline
\vertical W(H_3) & 120     & 2 & 3 & 4 & 10, 6, 1 & 3 & 2 \\
\vertical        & 2^3.3.5 & 3 & 2 & 2 & 10, 2    & 2 & 1 \\
\vertical        &         & 5 & 3 & 3 & 10, 4, 2 & 3 & 1 \\
\hline
\vertical W(H_4) & 14400       & 2 & 4 & 7 & 34, 24, 9, 1     & 4 & 3 \\
\vertical        & 2^6.3^2.5^2 & 3 & 3 & 3 & 34, 11, 2        & 3 & 1 \\
\vertical        &             & 5 & 5 & 6 & 34, 20, 11, 4, 2 & 5 & 2 \\
\hline
\end{array}$$

\newpage

$$\begin{array}{|c|c|@{\!}c@{\!}||c|c|l|c|c|c|}
\hline
\vphantom{\DS{\frac{A}{B}}} 
G & |G| & \quad p\quad & \loewy_p(G) & S_p(G) & d_0, d_1, d_2, \dots & 
\loewy_p(G,1) & \ext_p^1(G,1)\\
\hline
\hline
\vertical \AG_5 & 60      & 2 & 2 & 2 & 5, 1    & 2 & 1 \\
\vertical       & 2^2.3.5 & 3 & 2 & 2 & 5, 1    & 2 & 1 \\
\vertical       &         & 5 & 3 & 3 & 5, 2, 1 & 3 & 1 \\ 
\hline 
\vertical \AG_6 & 360       & 2 & 3 & 3 & 7, 2, 1 & 3 & 1 \\
\vertical       & 2^3.3^2.5 & 3 & 3 & 3 & 7, 2, 1 & 3 & 1 \\
\vertical       &           & 5 & 3 & 3 & 7, 2, 1 & 3 & 1 \\
\hline 
\vertical \AG_7 & 2520        & 2 & 3 & 3 & 9, 3, 1 & 3 & 1 \\
\vertical       & 2^3.3^2.5.7 & 3 & 3 & 3 & 9, 3, 1 & 3 & 1 \\
\vertical       &             & 5 & 2 & 2 & 9, 1    & 2 & 1 \\
\vertical       &             & 7 & 3 & 3 & 9, 2, 1 & 3 & 1 \\
\hline
\vertical \AG_8 & 20160       & 2 & 4 & 5 & 14, 6, 2, 1 & 4 & 2 \\
\vertical       & 2^6.3^2.5.7 & 3 & 3 & 3 & 14, 6, 2    & 3 & 1 \\
\vertical       &             & 5 & 3 & 3 & 14, 3, 1    & 2 & 1 \\
\vertical       &             & 7 & 3 & 3 & 14, 2, 1    & 3 & 1 \\
\hline 
\vertical \AG_9 & 181440      & 2 & 4 & 5 & 18, 8, 3, 1  & 4 & 2 \\
\vertical       & 2^6.3^4.5.7 & 3 & 4 & 6 & 18, 10, 3, 1 & 4 & 3 \\
\vertical       &             & 5 & 3 & 3 & 18, 4, 1     & 2 & 1 \\
\vertical       &             & 7 & 2 & 2 & 18, 1        & 2 & 1 \\
\hline 
\vertical \AG_{10} & 1814400       & 2 & 5 & 7 & 24, 12, 6, 2, 1 & 5 & 2 \\
\vertical          & 2^7.3^4.5^2.7 & 3 & 4 & 6 & 24, 13, 4, 1    & 4 & 3 \\
\vertical          &               & 5 & 3 & 3 & 24, 4, 1        & 3 & 1 \\
\vertical          &               & 7 & 3 & 3 & 24, 3, 1        & 2 & 1 \\
\hline 
\vertical \AG_{11} & 19958400         & 2 & 5 & 7 & 31, 17, 8, 3, 1 & 5 & 2 \\
\vertical          & 2^7.3^4.5^2.7.11 & 3 & 4 & 5 & 31, 16, 6, 1    & 4 & 2 \\
\vertical          &                  & 5 & 3 & 3 & 31, 6, 1        & 3 & 1 \\
\vertical          &                  & 7 & 3 & 3 & 31, 4, 1        & 2 & 1 \\
\vertical          &                  & 11& 3 & 3 & 31, 2, 1        & 3 & 1 \\
\hline
\vertical \AG_{12} & 239500800        & 2 & 6 & 10 & 43, 25, 13, 6, 2, 1 & 6 & 2 \\
\vertical          & 2^9.3^5.5^2.7.11 & 3 & 5 & 8 & 43, 22, 9, 2, 1      & 5 & 3 \\
\vertical          &                  & 5 & 3 & 3 & 43, 10, 2            & 3 & 1 \\
\vertical          &                  & 7 & 3 & 3 & 43, 5, 1             & 2 & 1 \\
\vertical          &                  & 11& 3 & 3 & 43, 2, 1             & 3 & 1 \\
\hline
\end{array}$$

%

\newpage
{\small 
$$\begin{array}{|c|c|@{\!}c@{\!}||c|c|l|c|c|}
\hline
\vphantom{\DS{\frac{A}{B}}} 
G & |G| & \quad p\quad & \loewy_p(G) & S_p(G) & d_0, d_1, d_2, \dots & 
\loewy_p(G,1) & \ext_p^1(G,1) \\
\hline
\hline
\vertical GL(3,2) & 168     & 2 & 3 & 3 & 6, 2, 1 & 3 & 1\\
\vertical         & 2^3.3.7 & 3 & 2 & 2 & 6, 1    & 2 & 1\\
\vertical         &         & 7 & 3 & 3 & 6, 2, 1 & 3 & 1\\
\hline
\vertical SL(2,8) & 504       & 2 & 2 & 2 & 9, 1          & 2 & 1\\
\vertical         & 2^3.3^2.7 & 3 & 5 & 5 & 9, 4, 3, 2, 1 & 5 & 1\\
\vertical         &           & 7 & 4 & 4 & 9, 3, 2, 1    & 4 & 1\\
\hline
\vertical SL(3,3) & 5616       & 2 & 5 & 5 & 12, 5, 3, 2, 1 & 5 & 1\\
\vertical         & 2^4.3^3.13 & 3 & 3 & 3 & 12, 3, 1       & 3 & 1\\
\vertical         &            & 13& 5 & 5 & 12, 4, 3, 2, 1 & 5 & 1\\
\hline 
\vertical SU(3,3) & 6048      & 2 & 6 & 7 & 14, 9, 6, 4, 2, 1 & 6 & 2\\
\vertical         & 2^5.3^3.7 & 3 & 3 & 3 & 14, 5, 1          & 3 & 1\\
\vertical         &           & 7 & 3 & 3 & 14, 2, 1          & 3 & 1\\
\hline
\vertical M_{11} & 7920         & 2 & 5 & 5 & 10, 5, 3, 2, 1 & 5 & 1\\
\vertical        & 2^4.3^2.5.11 & 3 & 2 & 2 & 10, 2          & 2 & 1\\
\vertical        &              & 5 & 2 & 2 & 10, 1          & 2 & 1\\
\vertical        &              & 11& 3 & 3 & 10, 2, 1       & 3 & 1\\
\hline
\vertical P\!Sp(4,3) & 25920     & 2 & 4 & 5 & 20, 12, 5, 1    & 4 & 2\\
\vertical          & 2^6.3^4.5 & 3 & 5 & 7 & 20, 14, 8, 3, 1 & 5 & 2\\
\vertical          &           & 5 & 2 & 2 & 20, 1           & 2 & 1\\
\hline
\vertical M_{12} & 95040        & 2 & 4 & 7 & 15, 9, 3, 1 & 4 & 3\\
\vertical        & 2^6.3^3.5.11 & 3 & 3 & 3 & 15, 4, 1    & 3 & 1\\
\vertical        &              & 5 & 2 & 2 & 15, 2       & 2 & 1\\
\vertical        &              & 11& 3 & 3 & 15, 2, 1   & 3 & 1\\
\hline 
\vertical J_1 & 175560          & 2 & 2 & 2 & 15, 4       & 2 & 1\\
\vertical     & 2^3.3.5.7.11.19 & 3 & 2 & 2 & 15, 4       & 2 & 1\\
\vertical     &                 & 5 & 3 & 3 & 15, 6, 3    & 3 & 1\\
\vertical     &                 & 7 & 2 & 2 & 15, 1       & 2 & 1\\
\vertical     &                 & 11& 2 & 2 & 15, 1       & 2 & 1\\
\vertical     &                 & 19& 4 & 4 & 15, 3, 2, 1 & 4 & 1\\
\hline 
\vertical M_{22} & 443520         & 2 & 4 & 5 & 12, 5, 2, 1 & 4 & 2\\
\vertical        & 2^7.3^2.5.7.11 & 3 & 2 & 2 & 12, 2       & 2 & 1\\
\vertical        &                & 5 & 2 & 2 & 12, 1       & 2 & 1\\
\vertical        &                & 7 & 3 & 3 & 12, 2, 1    & 3 & 1\\
\vertical        &                & 11& 3 & 3 & 12, 2, 1    & 3 & 1\\
\hline
\vertical J_2 & 604800        & 2 & 4 & 5 & 21, 11, 3, 1 & 4 & 2\\
\vertical     & 2^7.3^3.5^2.7 & 3 & 3 & 3 & 21, 7, 1,       & 3 & 1\\
\vertical     &               & 5 & 5 & 5 & 21, 10, 6, 2, 1 & 5 & 1\\
\vertical     &               & 7 & 2 & 2 & 21, 1           & 2 & 1\\
\hline 
\vertical H\!S & 44352000         & 2 & 5 & 9 & 24, 15, 8, 3, 1 & 5 & 3\\
\vertical    & 2^9.3^2.5^3.7.11 & 3 & 2 & 2 & 24, 5 & 2 & 1\\
\vertical    &                  & 5 & 3 & 4 & 24, 8, 2 & 3 & 2\\
\vertical    &                  & 7 & 2 & 2 & 24, 1 & 2 & 1\\
\vertical    &                  & 11 & 3 & 3 & 24, 2, 1 & 3 & 1\\
\hline
\end{array}$$

%

\newpage

$$\begin{array}{|c|c|@{\!}c@{\!}||c|c|l|c|c|}
\hline
\vphantom{{\frac{A}{B}}} 
G & |G| & \quad p\quad & \loewy_p(G) & S_p(G) & d_0, d_1, d_2, \dots & 
\loewy_p(G,1) & \ext_p^1(G,1) \\
\hline
\hline
\vertical P\!S\!L(2,2) & 6   & 2 & 2 & 2 & 3, 1 & 2       & 1 \\
\vertical \simeq\SG_3        & 2.3 & 3 & 2 & 2 & 3, 1 & 2 & 1 \\
\hline
\vertical P\!S\!L(2,3) & 12    & 2 & 2 & 2 & 4, 1    & 2 & 1 \\
\vertical \simeq\AG_4        & 2^2.3 & 3 & 3 & 3 & 4, 2, 1 & 3 & 1 \\
\hline
\vertical P\!S\!L(2,4) & 60      & 2 & 2 & 2 & 5, 1    & 2 & 1 \\
\vertical \simeq P\!S\!L(2,5) & 2^2.3.5 & 3 & 2 & 2 & 5, 1    & 2 & 1 \\
\vertical \simeq\AG_5        &         & 5 & 3 & 3 & 5, 2, 1 & 3 & 1 \\
\hline 
\vertical P\!S\!L(2,7) & 168     & 2 & 3 & 3 & 6, 2, 1 & 3 & 1 \\
\vertical              & 2^3.3.7 & 3 & 2 & 2 & 6, 1    & 2 & 1 \\
\vertical              &         & 7 & 3 & 3 & 6, 2, 1 & 3 & 1 \\
\hline
\vertical P\!S\!L(2,8) & 504       & 2 & 2 & 2 & 9, 1          & 2 & 1 \\
\vertical              & 2^3.3^2.7 & 3 & 5 & 5 & 9, 4, 3, 2, 1 & 5 & 1 \\
\vertical              &           & 7 & 4 & 4 & 9, 3, 2, 1    & 4 & 1 \\
\hline
\vertical P\!S\!L(2,9) & 360       & 2 & 3 & 3 & 7, 2, 1 & 3 & 1 \\
\vertical \simeq \AG_6        & 2^3.3^2.5 & 3 & 3 & 3 & 7, 2, 1 & 3 & 1 \\
\vertical              &           & 5 & 3 & 3 & 7, 2, 1 & 3 & 1 \\
\hline
\vertical P\!S\!L(2,11) & 660        & 2 & 2 & 2 & 8, 2    & 2 & 1 \\
\vertical               & 2^2.3.5.11 & 3 & 2 & 2 & 8, 2    & 2 & 1 \\
\vertical               &            & 5 & 3 & 3 & 8, 2, 1 & 3 & 1 \\
\vertical               &            & 11& 3 & 3 & 8, 2, 1 & 3 & 1 \\
\hline
\vertical P\!S\!L(2,13) & 1092       & 2 & 2 & 2 & 9, 2       & 2 & 1 \\
\vertical               & 2^2.3.7.13 & 3 & 2 & 2 & 9, 2       & 2 & 1 \\
\vertical               &            & 7 & 4 & 4 & 9, 3, 2, 1 & 4 & 1 \\
\vertical               &            & 13& 3 & 3 & 9, 2, 1    & 3 & 1 \\
\hline
\vertical P\!S\!L(2,16) & 4080       & 2 & 2 & 2 & 17, 1                      & 2 & 1 \\
\vertical               & 2^4.3.5.17 & 3 & 3 & 3 & 17, 5, 2                   & 2 & 1 \\
\vertical               &            & 5 & 5 & 5 & 17, 6, 4, 2, 1             & 3 & 1 \\
\vertical               &            & 17& 9 & 9 & 17, 8, 7, 6, 5, 4, 3, 2, 1 & 9 & 1 \\
\hline
\vertical P\!S\!L(2,17) & 2448       & 2 & 5 & 5 & 11, 4, 3, 2, 1 & 5 & 1 \\
\vertical               & 2^4.3^2.17 & 3 & 5 & 5 & 11, 4, 3, 2, 1 & 5 & 1 \\
\vertical               &            & 17& 3 & 3 & 11, 2, 1       & 3 & 1 \\
\hline
\vertical P\!S\!L(2,19) & 3420         & 2 & 2 & 2 & 12, 3          & 2 & 1 \\
\vertical               & 2^2.3^2.5.19 & 3 & 5 & 5 & 12, 4, 3, 2, 1 & 5 & 1 \\
\vertical               &              & 5 & 3 & 3 & 12, 4, 2       & 3 & 1 \\
\vertical               &              & 19& 3 & 3 & 12, 2, 1       & 3 & 1 \\
\hline
\vertical P\!S\!L(2,23) & 6072        & 2 & 4 & 4 & 14, 5, 3, 1       & 3 & 1 \\
\vertical               & 2^3.3.11.23 & 3 & 3 & 3 & 14, 4, 1          & 2 & 1 \\
\vertical               &             & 11& 6 & 6 & 14, 5, 4, 3, 2, 1 & 6 & 1 \\
\vertical               &             & 23& 3 & 3 & 14, 2, 1          & 3 & 1 \\
\hline 
\vertical P\!S\!L(2,25) & 7800         & 2 & 4 & 4 & 15, 5, 3, 1          & 3 & 1 \\
\vertical               & 2^3.3.5^2.13 & 3 & 3 & 3 & 15, 4, 1             & 2 & 1 \\
\vertical               &              & 5 & 3 & 3 & 15, 2, 1             & 3 & 1 \\
\vertical               &              & 13& 7 & 7 & 15, 6, 5, 4, 3, 2, 1 & 7 & 1 \\
\hline
\vertical P\!S\!L(2,27) & 9828         & 2 & 2 & 2 & 16, 4                & 2 & 1 \\
\vertical               & 2^2.3^3.7.13 & 3 & 3 & 3 & 16, 2, 1             & 3 & 1 \\
\vertical               &              & 7 & 4 & 4 & 16, 6, 4, 2          & 4 & 1 \\
\vertical               &              & 13& 7 & 7 & 16, 6, 5, 4, 3, 2, 1 & 7 & 1 \\
\hline     
\end{array}$$}





\begin{thebibliography}{DLMM}

\bibitem[B]{bonnafe sym} {\sc C. Bonnaf\'e}, A note on the Grothendieck ring 
of the symmetric group, to appear in {\it  C. R. Math. Acad. Sci. Paris}.

\bibitem[CR]{curtis} {\sc C.W. Curtis \& I. Reiner}, 
{\it Methods of representation thoery, vol. I}, Wiley, New York, 1981; 
reprinted 1990 as Wiley Classics Library.

\bibitem[DM]{dimibook} {\sc F. Digne \& J. Michel}, 
{\it Representations of finite groups of Lie type}, 
London math. soc. student texts {\bf 21}, 
Cambridge university press, 1991.

\bibitem[GAP3]{GAP} {\sc Martin Sch\"onert et. al.}, 
GAP -- Groups, Algorithms, and Programming -- version 3 release 4 patchlevel 4. 
Lehrstuhl D f\"ur Mathematik, Rheinisch Westf\"alische Technische Hochschule, 
Aachen, Germany, 1997.

\bibitem[GP]{geck} {\sc M. Geck \& G. Pfeiffer}, {\it 
Characters of finite Coxeter groups and Iwahori-Hecke algebras}, 
London Math. Soc. Monographs, New series no. {\bf 21}, 
Oxford University Press, 2000.

\bibitem[H]{higman} {\sc G. Higman}, The units of group-rings, 
{\it Proc. London Math. Soc.} {\bf 46} (1940), 231--248.

\bibitem[I]{isaacs} {\sc M.I. Isaacs}, {\it  Character theory of finite groups}, 
Pure and Applied Mathematics {\bf 69}, 
Academic Press, New York-London, 1976.

\bibitem[L]{lusztig} {\sc G. Lusztig}, Leading coefficients of character values 
of Hecke algebras, {\it The Arcata Conference on Representations of Finite Groups 
(Arcata, Calif., 1986)}, 235--262, Proc. Sympos. Pure Math. {\bf 47}, Part 2, 
Amer. Math. Soc., Providence, RI, 1987. 

\bibitem[S]{saksonov} {\sc A.I. Saksonov}, 
The integral ring of characters of a finite group (Russian), 
Vesci Akad. Navuk BSSR Ser. Fiz.-Mat. Navuk {\bf 1966} (1966) no. 3, 69--76.
\end{thebibliography}
\end{document}